\documentstyle{amltd2004}
\begin{document}
\annalsline{158}{2003}
\received{June 11, 2001}
\startingpage{635}
\def\bye{\end{document}}
 \font\tenrm=cmr10
\def\ritem#1{\item[{\rm #1}]}
\catcode`\@=11
\font\twelvemsb=msbm10 scaled 1100
\font\tenmsb=msbm10
\font\ninemsb=msbm10 scaled 800
\newfam\msbfam
\textfont\msbfam=\twelvemsb  \scriptfont\msbfam=\ninemsb
  \scriptscriptfont\msbfam=\ninemsb
\def\msb@{\hexnumber@\msbfam}
\def\Bbb{\relax\ifmmode\let\next\Bbb@\else
 \def\next{\errmessage{Use \string\Bbb\space only in math
mode}}\fi\next}
\def\Bbb@#1{{\Bbb@@{#1}}}
\def\Bbb@@#1{\fam\msbfam#1}
\catcode`\@=12

 \catcode`\@=11
\font\twelveeuf=eufm10 scaled 1100
\font\teneuf=eufm10
\font\nineeuf=eufm7 scaled 1100
\newfam\euffam
\textfont\euffam=\twelveeuf  \scriptfont\euffam=\teneuf
  \scriptscriptfont\euffam=\nineeuf
\def\euf@{\hexnumber@\euffam}
\def\frak{\relax\ifmmode\let\next\frak@\else
 \def\next{\errmessage{Use \string\frak\space only in math
mode}}\fi\next}
\def\frak@#1{{\frak@@{#1}}}
\def\frak@@#1{\fam\euffam#1}
\catcode`\@=12

\newcommand{\ovl}{\overline}
\newcommand{\n}{\noindent}
\newcommand{\vp}{\varepsilon}
\newcommand{\cl}[1]{{{\cal{#1}}}}
\newcommand{\bb}[1]{{\Bbb{ #1}}} 
\newcommand{\tr}{{\rm  tr}}
\newcommand{\intl}{\int\limits}



 \font\eighteensyb=cmr10 scaled \magstep3
\title{Hochschild cohomology of factors\\ with property  \hbox{\eighteensyb \char0}}
\shorttitle{Hochschild cohomology} 

\newif\ifauthor
\def\author#1{\vskip15pt
\hbox to\hsize{\hss\tenrm By \tensc#1\ifacks\global\acksfalse*\fi\hss}
\ifshort\else\xdef\theauthors{{\eightsc\uppercase{E. Christensen, F. Pop, A. M. Sinclair,  and  R. R. Smith}}}\fi%
\vskip7pt
\vskip\baselineskip
\global\authortrue\everypar={\global\authorfalse\everypar={}}}

 \author{Erik Christensen$^\ast$, Florin Pop, Allan M.~Sinclair,}

\vglue-20pt \centerline{{\tenrm and} {\tensc Roger R.~Smith$^{\dagger}$}}

 \institutions{Institute for Mathematiske Fag, University of Copenhagen, Copenhagen,
Denmark\\
{\eightpoint {\it E-mail address\/}:
  echris@math.ku.dk}\\
\vglue6pt
Wagner College, Staten Island, NY\\
{\eightpoint {\it E-mail address\/}:
 fpop@wagner.edu}\\
\vglue6pt
University of Edinburgh, Edinburgh, Scotland, United Kingdom\\
{\eightpoint {\it E-mail address\/}:
 allan@maths.ed.ac.uk}\\
\vglue6pt
Texas A\&M University, College Station, TX \\
{\eightpoint {\it E-mail address\/}: rsmith@math.tamu.edu}}
 
\begin{center}{{\it Dedicated to the memory of Barry Johnson}, 1937--2002}

\end{center}

\vglue8pt \centerline{\bf Abstract}
\vglue8pt
The main result of this paper is that the $k^{\rm th}$ continuous Hochschild
cohomology groups $H^k(\cl M,\cl M)$ and $H^k(\cl M,B(H))$ of a von~Neumann
factor ${\cl M}\subseteq B(H)$ of type ${\rm II}_1$ with property $\Gamma$ are zero
for all positive integers $k$. The method of proof involves the construction of
hyperfinite subfactors with special properties and a new inequality of
Grothendieck type for multilinear maps. We prove joint continuity in the
$\|\cdot\|_2$-norm of separately ultraweakly continuous multilinear maps, and
combine these results to reduce to the case of completely bounded cohomology
which is already solved.

\vglue18pt
\section{Introduction}\label{sec1}
\vglue-6pt

The continuous Hochschild cohomology of von~Neumann algebras was initiated by
Johnson, Kadison and Ringrose in a series of papers \cite{JKR}, \cite{KR1}, \cite{KR2}
where they developed the basic theorems and techniques of the subject. 
>From their results, and from those of subsequent authors, it was natural 
to conjecture
that the $k^{\rm th}$ continuous Hochschild cohomology group $H^k(\cl M,\cl
M)$ of a von~Neumann algebra over itself is zero for all positive integers
$k$. This was verified by Johnson, Kadison and Ringrose, \cite{JKR}, for all
hyperfinite von~Neumann algebras and the cohomology was shown to split over
the center. A technical version of their result has been used in all
subsequent proofs and is applied below. Triviality of the cohomology groups
has interesting structural implications for von~Neumann algebras,
\cite[Chapter 7]{SS1} (which surveys the original work in this area by
Johnson, \cite{J}, and Raeburn and Taylor, \cite{RT}), and so it is important
to determine when this occurs.
\vfill
\footnoterule
 {\footnotesize  $^\ast$ Partially supported by a Scheme 4
collaborative grant from the London Mathematical Society.}
\vglue-2pt

 {\footnotesize $^{\dagger}$ Partially supported by a grant from the National
Science Foundation.} 
\eject

The representation theorem for completely bounded multilinear maps,
\cite{CS1}, which expresses such a map as a product of $*$-homomorphisms and
interlacing operators, was used by the first and third authors to show that
the completely bounded cohomology $H^k_{\rm cb}(\cl M,\cl M)$ is always zero
\cite{CS3}, \cite{CS4}, \cite{SS1}. Subsequently it was observed in \cite{SS2}, \cite{SS3}, \cite{SS4} that to
show that $H^k(\cl M, \cl M) = 0$, it suffices to reduce a normal cocycle to a
cohomologous one that is completely bounded in the first or last variable
only, while holding fixed the others. The multilinear maps that are completely
bounded in the first (or last) variable do not form a Hochschild complex;
however it is easier to check complete boundedness in one variable only
\cite{SS2}. In joint work with Effros, \cite{CES}, the first and third authors
had shown that if the type ${\rm II}_1$ central summand of a von~Neumann algebra
$\cl M$ is stable under tensoring with the hyperfinite type ${\rm II}_1$ factor $\cl
R$, then
\begin{equation}\label{eq1.1}
H^k(\cl M,\cl M) = H^k_{\rm cb}(\cl M,\cl M) = 0,\qquad k\ge 2.
\end{equation}
This reduced the conjecture to type ${\rm II}_1$ von~Neumann algebras, and a further
reduction to those von~Neumann algebras with separable preduals was
accomplished in \cite[\S 6.5]{SS1}. We note that we restrict to $k\ge 2$,
since the case $k=1$, in a different formulation, is the question of whether
every derivation of a von~Neumann algebra into itself is inner, and this was
solved independently by Kadison and Sakai, \cite{K}, \cite{S}.

The noncommutative Grothendieck inequality for normal bilinear forms on a
von~Neumann algebra due to Haagerup, \cite{H}  (but building on earlier work
of Pisier, \cite{Pi}) and the existence of hyperfinite subfactors with trivial
relative commutant due to Popa, \cite{Po}, have been the main tools for 
showing that suitable cocycles are completely bounded in the first variable,
\cite{C3}, \cite{SS2}, \cite{SS3}, \cite{SS4}. The importance of this inequality for derivation
problems on von~Neumann and $C^*$-algebras was initially observed in the work
of Ringrose, \cite{R}, and of the first author, \cite{C1}. The current state
of knowledge for the cohomology conjecture for type ${\rm II}_1$ factors may be
summarized as follows: \vglue-22pt
\phantom{bugger}
\begin{itemize}
\item[(i)] $\cl M$ is stable under tensoring by the hyperfinite type ${\rm II}_1$
factor $\cl R$, $k\ge 2$, \cite{CES};
\item[(ii)] $\cl M$ has property $\Gamma$ and $k=2$, \cite{C3}, \cite{CS3};
\item[(iii)] $\cl M$ has a Cartan subalgebra, \cite[$k=2$]{PS},
\cite[$k=3$]{CPSS}, \cite[$k\ge 2$]{SS2, SS3};
\item[(iv)] $\cl M$ has various technical properties relating to its action on
$L^2(\cl M, \tr)$ for $k=2$, (\cite{PS}), and conditions of this type were
verified for various classes of factors by Ge and Popa,  \cite{GP}.
\end{itemize}
\phantom{bugger}
\vglue-22pt
 
 The two test questions for the type ${\rm II}_1$ factor case are the following. Is
$H^k(\cl M, \cl M)$ equal to zero for factors with property $\Gamma$, and is
$H^2(VN(\bb F_2)$, $VN(\bb F_2))$ equal to zero for the von~Neumann factor of
the free group on two generators? The second is still open at this time; the
purpose of this paper is to give a positive answer to the first
(Theorems~\ref{thm6.2} and \ref{thm7.2}). If we change the coefficient module
to be any containing $B(H)$, then the question arises of whether analogous
results for $H^k({\cl M},B(H))$ are valid (see \cite{CES}). We will see below
that our methods are also effective in this latter case. 

The algebras of (i) above are called McDuff factors, since they were
studied in \cite{Mc1}, \cite{Mc2}. The hyperfinite factor $\mathcal R$ satisfies property
$\Gamma$ (defined in the next section), and it is an easy consequence of
the definition that the tensor product of an arbitrary type ${\rm II}_1$ factor
with a $\Gamma$-factor also has property $\Gamma$. Thus, as is well
known, the McDuff factors all have property $\Gamma$, and so the results
of this paper recapture the vanishing of cohomology for this class, \cite{CES}.
However, as was shown by Connes, \cite{Co}, the class of factors with property
$\Gamma$ is much wider. This was confirmed in recent work of Popa, \cite{Po2},
who constructed a family of $\Gamma$-factors with trivial fundamental
group. This precludes the possibility that they are McDuff factors, all of which
have fundamental group equal to ${\Bbb R}^+$.

The most general class of type ${\rm II}_1$ factors for which vanishing of
cohomology has been obtained is described in (iii). While there is some
overlap between those factors with Cartan subalgebras and those with
property $\Gamma$, the two classes do not appear to be directly related,
since their definitions are quite different.
It is not difficult to
verify that the infinite tensor product of an arbitrary sequence of type
${\rm II}_1$ factors has property $\Gamma$, using the $\|\cdot\|_2$-norm
density of the span of elements of the form
$x_1\otimes x_2\otimes \cdots \otimes x_n \otimes 1 \otimes 1\cdots$.
Voiculescu, \cite{V}, has exhibited a family of factors (which includes
$VN({\Bbb F}_2)$) having no Cartan subalgebras, but also failing to have
property $\Gamma$. This suggests that the infinite tensor product of
copies of this algebra might be an example of a factor with property
$\Gamma$ but without a Cartan subalgebra. This is unproved, and indeed the
question of whether
 $VN({\Bbb F}_2) \overline{\otimes} VN({\Bbb F}_2)$
has a Cartan subalgebra appears to be open at this time. While we do not
know of a factor with property $\Gamma$ but with no Cartan subalgebra,
these remarks indicate that such an example may well exist. Thus the
results of this paper and the earlier results of \cite{SS2} should be viewed as
complementary to one another, but not necessarily linked.

We now give a brief description of our approach to this problem; definitions
and a more extensive discussion of background material will follow in the next
section. For a factor $\cl M$ with separable predual and property $\Gamma$, we
construct a hyperfinite subfactor $\cl R\subseteq \cl M$ with trivial relative
commutant which enjoys the additional property of containing an asymptotically
commuting family of projections for the algebra $\cl M$ (fifth section). In
the third section we prove a Grothendieck inequality for $\cl R$-multimodular
normal multilinear maps, and in the succeeding section we show that separate
normality leads to joint continuity in the $\|\cdot\|_2$-norm (or,
equivalently, joint ultrastrong$^*$ continuity) on the closed unit ball of
$\cl M$. These three results are sufficient to obtain vanishing cohomology for
the case of a separable predual (sixth section), and we give the extension to
the general case at the end of the paper.

We refer the reader to our lecture notes on cohomology, \cite{SS1}, for many
of the results used here and to \cite{C2}, \cite{Co}, \cite{D}, \cite{Mc1}, \cite{Mc2}, \cite{MvN} for other material
concerning property $\Gamma$. We also take the opportunity to thank Professors
I.~Namioka and Z.~Piotrowski for their guidance on issues related to the
fourth section of the paper.  
\section{Preliminaries}\label{sec2}

Throughout the paper $\cl M$ will denote a type ${\rm II}_1$ factor with unique
normalized normal trace ${\rm tr}$. We write $\|x\|$ for the operator norm of an
element $x\in\cl M$, and $\|x\|_2$ for the quantity $({\rm tr}(xx^*))^{1/2}$, which
is the norm induced by the inner product $\langle x,y\rangle = {\rm tr}(y^*x)$ on
$\cl M$.

Property $\Gamma$ for a type ${\rm II}_1$ factor $\cl M$ was introduced by Murray
and\break von Neumann, \cite{MvN}, and is defined by the following requirement: \
given $x_1,\ldots,\break x_m\in \cl M$ and $\vp>0$, there exists a unitary $u\in
\cl M$, ${\rm tr}(u) = 0$, such that
\begin{equation}\label{eq2.1}
\|ux_j-x_ju\|_2 < \vp,\qquad 1\le j\le m.
\end{equation}
Subsequently we will use both this definition and the following equivalent
formulation due to Dixmier, \cite{D}. Given $\vp>0$, elements $x_1,\ldots,
x_m\in \cl M$, and a positive integer $n$, there exist orthogonal projections
$\{p_i\}^n_{i=1}\in \cl M$, each of trace $n^{-1}$ and summing to 1, such that
\begin{equation}\label{eq2.2}
\|p_ix_j-x_jp_i\|_2 < \vp,\qquad 1\le j \le m,\quad 1\le i\le n.
\end{equation}
In \cite{Po}, Popa showed that each type ${\rm II}_1$ factor $\cl M$ with separable\break
predual contains a hyperfinite subfactor $\cl R$ with trivial relative
commutant\break  $(\cl R'\cap \cl M = \bb C 1)$, answering positively an earlier
question posed by Kadison. In the presence of property $\Gamma$, we will
extend Popa's theorem by showing that  $\cl R$ may be chosen to contain,
within a maximal abelian subalgebra, projections which satisfy (\ref{eq2.2}).
This result is Theorem~\ref{thm5.3}.

We now briefly recall the definition of continuous Hochschild cohomology for
von~Neumann algebras. Let $\cl X$ be a Banach $\cl M$-bimodule and let $\cl
L^k(\cl M,\cl X)$ be the Banach space of $k$-linear bounded maps from the
$k$-fold Cartesian product $\cl M^k$ into $\cl X$, $k\ge 1$. For $k=0$, we
define $\cl L^0(\cl M,\cl X)$ to be $\cl X$. The coboundary operator
$\partial^k\colon \ \cl L^k(\cl M,\cl X)\to \cl L^{k+1}(\cl M,\cl X)$ (usually
abbreviated to just $\partial$) is defined, for $k\ge 1$, by
\begin{eqnarray}
\qquad\quad\partial\phi(x_1,\ldots, x_{k+1}) &= &x_1\phi(x_2,\ldots, x_{k+1})\label{eq2.3}
\\
&& +\ \sum^k_{i=1} (-1)^i \phi(x_1,\ldots, x_{i-1}, x_ix_{i+1},
x_{i+2},\ldots, x_{k+1})\nonumber\\
&& +\ (-1)^{k+1} \phi(x_1,\ldots, x_k)x_{k+1}, \nonumber
\end{eqnarray}
for $x_1,\ldots, x_{k+1}\in \cl M$. When $k=0$, we define $\partial\xi$, for
$\xi\in \cl X$, by
\begin{equation}\label{eq2.4}
\partial\xi(x) = x\xi -\xi x,\qquad x\in\cl M.
\end{equation}
It is routine to check that $\partial^{k+1} \partial^k = 0$, and so $\hbox{Im
} \partial^k$ (the space of coboundaries) is contained in $\hbox{Ker }
\partial^{k+1}$ (the space of cocycles). The continuous Hochschild cohomology
groups $H^k(\cl M,\cl X)$ are then defined to be the quotient vector spaces
$\hbox{Ker } \partial^k/\hbox{Im } \partial^{k-1}$, $k\ge 1$. When $\cl X$ is
taken to be $\cl M$, an element $\phi\in \cl L^k(\cl M,\cl M)$ is normal if
$\phi$ is separately continuous in each of its variables when both range and
domain are endowed with the ultraweak topology induced by the predual~$\cl
M_*$.

Let $\cl N\subseteq \cl M$ be a von~Neumann subalgebra, and assume that $\cl
M$ is represented on a Hilbert space $H$. Then $\phi\colon \ \cl M^k\to B(H)$
is $\cl N$-multimodular if the following conditions are satisfied by all $a\in
\cl N$, $x_1,\ldots, x_k\in\cl M$, and $1\le i \le k-1$:
\begin{eqnarray}
\label{eq2.5}
a\phi(x_1,\ldots, x_k) &= &\phi(ax_1,x_2,\ldots, x_k),\\
\label{eq2.6}
\phi(x_1,\ldots, x_k)a &= &\phi(x_1,\ldots, x_{k-1}, x_ka),\\
\label{eq2.7}
\phi(x_1,\ldots,  x_ia, x_{i+1},\ldots, x_k) &=& \phi(x_1,\ldots, x_i,
ax_{i+1},\ldots, x_k).
\end{eqnarray}
A fundamental result of Johnson, Kadison and Ringrose, \cite{JKR}, states
that each cocycle $\phi$ on $\cl M$ is cohomologous to a normal cocycle $\phi
- \partial\psi$, which can also be chosen to be $\cl N$-multimodular for any
given hyperfinite subalgebra $\cl N\subseteq \cl M$. This has been the
starting point for all subsequent theorems in von~Neumann algebra cohomology,
since it permits the substantial simplification of considering only $\cl
N$-multimodular normal cocycles for a suitably chosen hyperfinite subalgebra
$\cl N$, \cite[Chapter 3]{SS1}. The present paper will provide another
instance of this.

The matrix algebras $\bb M_n(\cl M)$ over a von~Neumann algebra (or
$C^*$-algebra) $\cl M$ carry natural $C^*$-norms inherited from $\bb
M_n(B(H)) = B(H^n)$, when $\cl M$ is represented on $H$. Each bounded map
$\phi\colon \ \cl M\to B(H)$ induces a sequence of maps $\phi^{(n)}\colon \
\bb M_n(\cl M)\to \bb M_n(B(H))$ by applying $\phi$ to each matrix entry\break
(it is usual to denote these by $\phi_n$ but we have adopted $\phi^{(n)}$ to
avoid notational difficulties in the sixth section). Then $\phi$ is said to be
completely bounded if $\sup\limits_{n\ge 1} \|\phi^{(n)}\| < \infty$, and this
supremum defines the completely bounded norm $\|\phi\|_{\rm cb}^{\phantom{|}}$ (see \cite{ER}, \cite{Pa}
for the extensive theory of such maps). A parallel theory for multilinear maps
was developed in \cite{CS1}, \cite{CS2}, using $\phi\colon \ \cl M^k\to \cl M$ to
replace the product in matrix multiplication. We illustrate this with $k=2$.
The $n$-fold amplification $\phi^{(n)}\colon \ \bb M_n(\cl M)\times \bb
M_n(\cl M) \to \bb M_n(\cl M)$ of a bounded bilinear map $\phi\colon \ \cl
M\times \cl M\to \cl M$ is defined as follows. For matrices $(x_{ij})$,
$(y_{ij}) \in \bb M_n(\cl M)$, the $(i,j)$ entry of $\phi^{(n)}((x_{ij}),
(y_{ij}))$ is $\sum\limits^n_{k=1} \phi(x_{ik}, y_{kj})$. We note \pagebreak that if
$\phi$ is\break $\cl N$-multimodular, then it is easy to verify from the definition
of $\phi^{(n)}$ that this map is $\bb M_n(\cl N)$-multimodular for each $n\ge
1$, and this will be used in the next section.

As before, $\phi$ is said to be completely bounded if $\sup\limits_{n\ge 1}
\|\phi^{(n)}\| < \infty$. By requiring all cocycles and coboundaries to be
completely bounded, we may define the completely bounded Hochschild
cohomology groups $H^k_{\rm cb}(\cl M,\cl M)$ and $H^k_{\rm cb}({\cl M},B(H)
)$ analogously to the continuous case. It was shown in \cite{CS3}, \cite{CS4} (see
also \cite[Chapter 4]{SS1}) that $H^k_{\rm cb}(\cl M,\cl M) = 0$ for $k\ge 1$ and
all von~Neumann algebras $\cl M$, exploiting the representation theorem for
completely bounded multilinear maps, \cite{CS1}, which is lacking in the
bounded case. This built on earlier work, \cite{CES}, on completely
bounded cohomology when the module is $B(H)$. Subsequent investigations have
focused on proving that cocycles are cohomologous to completely bounded ones,
\cite{CPSS}, \cite{PS}, or to ones which exhibit complete boundedness in one of the
variables \cite{C3}, \cite{SS2}, \cite{SS3}, \cite{SS4}. We will also employ this strategy here. 
 
\vglue-8pt
\section{A multilinear Grothendieck inequality}\label{sec3}
\vglue-6pt

The noncommutative Grothendieck inequality for bilinear forms, \cite{Pi}, and
its normal counterpart, \cite{H}, have played a fundamental role in Hochschild
cohomology theory \cite[Chapter 5]{SS1}. The main use has been to show that
suitable normal cocycles are completely bounded in at least one variable
\cite{CPSS}, \cite{SS2}, \cite{SS3}, \cite{SS4}. In this section we prove a multilinear version of
this inequality which will allow us to connect continuous and completely
bounded cohomology in the sixth section.

If $\cl M$ is a type ${\rm II}_1$ factor and $n$ is a positive integer, we denote by
$\tr_n$ the normalized trace on $\bb M_n(\cl M)$, and we introduce the
quantity $\rho_n(X) = (\|X\|^2 + n  \ tr_n(X^*X))^{1/2}$, for $X\in \bb
M_n(\cl M)$. We let $\{E_{ij}\}^n_{i,j=1}$ be the standard matrix units for
$\bb M_n$ ($\{e_{ij}\}^n_{i,j=1}$ is the more usual way of writing these
matrix units, but we have chosen upper case letters to conform to our
conventions on matrices). If $\phi^{(n)}$ is the $n$-fold amplification of the
$k$-linear map $\phi$ on $\cl M$ to $\bb M_n(\cl M)$, then
\vglue2pt
\centerline{$\phi^{(n)}(E_{11}X_1E_{11},\ldots, E_{11}X_kE_{11}),\qquad X_i\in \bb
M_n(\cl M),$}
\vglue4pt\noindent  is simply $\phi$ evaluated at the (1,1) entries of these
matrices, leading to the inequality \begin{equation}\label{eq3.1}
\|\phi^{(n)}(E_{11}X_1E_{11},\ldots, E_{11}X_kE_{11})\| \le \|\phi\| \|X_1\|
\ldots \|X_k\|.
\end{equation}
Our objective in Theorem~\ref{thm3.3} is to successively remove the matrix
units from (\ref{eq3.1}), moving from left to right, at the expense  of
increasing the right-hand side of this inequality. The following two variable
inequality will allow us to achieve this for certain multilinear maps.

\proclaim{Lemma}\label{lem3.1}
Let $\cl M\subseteq B(H)$ be a type ${\rm II}_1$ factor with a hyperfinite subfactor
$\cl N$ of trivial relative commutant{\rm ,} let $C>0$ and let $n$ be a positive
integer. If $\psi\colon \ \bb M_n(\cl M)\times \bb M_n(\cl M)\to B(H)$ is
a normal bilinear map satisfying
\begin{equation}\label{eq3.2}
\psi(XA,Y) = \psi(X,AY), \qquad A\in\bb M_n(\cl N),\qquad X,Y\in \bb M_n(\cl
M), \qquad\end{equation}
and
\begin{equation}\label{eq3.3}
\|\psi(XE_{11},E_{11}Y)\|\le C\|X\| \|Y\|,\qquad X,Y\in \bb M_n(\cl M),
\end{equation}
then
\begin{equation}\label{eq3.4}
\|\psi(X,Y)\| \le C\rho_n(X)\rho_n(Y),\qquad X,Y\in\bb M_n(\cl M).
\end{equation}
\endproclaim

\demo{Proof}
Let $\eta$ and $\nu$ be arbitrary unit vectors in $H^n$ and define a normal
bilinear form on $\bb M_n(\cl M)\times \bb M_n(\cl M)$ by
\begin{equation}\label{eq3.5}
\theta(X,Y) = \langle\psi(XE_{11}, E_{11}Y)\eta,\nu\rangle
\end{equation}
for $X,Y\in \bb M_n(\cl M)$. Then $\|\theta\| \le C$ by (\ref{eq3.3}). By the
noncommutative Grothendieck inequality for normal bilinear forms on a
von~Neumann algebra, \cite{H}, there exist normal states $f,F,g$ and $G$ on
$\bb M_n(\cl M)$ such that
\begin{equation}\label{eq3.6}
|\theta(X,Y)|\le C(f(XX^*) + F(X^*X))^{1/2} (g(YY^*) + G(Y^*Y))^{1/2} \hskip.4in
\end{equation}
for all $X,Y\in \bb M_n(\cl M)$. From (\ref{eq3.2}), (\ref{eq3.5}) and
(\ref{eq3.6}),
\begin{eqnarray}
|\langle\psi(X,Y)\eta,\nu\rangle| &=& \left|\sum^n_{j=1} \langle\psi (XE_{j1}
E_{11}, E_{11}E_{1j}Y)\eta,\nu\rangle\right|\label{eq3.7}
\\
&\le& \sum^n_{j=1} |\theta(XE_{j1}, E_{1j}Y)|,\nonumber
\end{eqnarray}
which we can then estimate by
\begin{eqnarray}\label{eq3.8}
&&\hskip-48pt  C\sum^n_{j=1} (f(XE_{j1} E_{1j}X^*) +  F(E_{1j}X^*XE_{j1}))^{1/2}\\
&&\qquad\quad\times\  (g(E_{1j}
YY^*E_{j1}) 
 + G(Y^*E_{j1}E_{1j}Y))^{1/2},\nonumber
\end{eqnarray} 
and this is at most
\begin{equation}\label{eq3.9}
C\left(f(XX^*)  +  \sum^n_{j=1} F(E_{1j}X^*XE_{j1})\right)^{1/2}
\left(\sum^n_{j=1} g(E_{1j}YY^*E_{j1}) + G(Y^*Y)\right)^{1/2},
\end{equation}
by the Cauchy-Schwarz \pagebreak inequality.

Let $\{\cl N_\lambda\}_{\lambda\in\Lambda}$ be an increasing net of matrix
subalgebras of $\cl N$ whose union is ultraweakly dense in $\cl N$. Let $\cl
U_\lambda$ denote the unitary group of $\bb M_n(\cl N_\lambda)$ with
normalized Haar measure $dU$. Since $\bb M_n(\cl N)' \cap \bb M_n(\cl M) = \bb
C1$, a standard argument (see \cite[5.4.4]{SS1}) gives
\begin{equation}\label{eq3.10}
\tr_n(X)1 = \lim_\lambda \int_{\cl U_\lambda} U^*XU\ dU
\end{equation}
in the ultraweak topology. Substituting $XU$ and $U^*Y$ respectively for $X$
and $Y$ in (\ref{eq3.7})--(\ref{eq3.9}), integrating over ${\cl U}_\lambda$ and
using the Cauchy-Schwarz inequality give
\begin{eqnarray}
&&\label{eq3.11}\\
|\langle\psi(X,Y)\eta,\nu\rangle| &=& |\langle\psi(XU,U^*Y)\eta,\nu\rangle|
\nonumber\\
&\le &C\left(f(XX^*) + \sum^n_{j=1} F\left(E_{1j} \int_{\cl U_\lambda}
U^*X^*XU \  dU \ E_{j1}\right)\right)^{1/2}
\nonumber \\
&&\times\ \left(\sum^n_{j=1} g\left(E_{1j}
\int_{\cl U_\lambda} U^*YY^*U \ dU \ E_{j1}\right) + G (Y^*Y)\right)^{1/2}.\nonumber
\end{eqnarray}
Now take the ultraweak limit over $\lambda\in\Lambda$ in (\ref{eq3.11}) to
obtain 
\begin{eqnarray}
 \qquad |\langle\psi(X,Y)\eta,\nu\rangle| 
& \le& C\left(f(XX^*) + \sum^n_{j=1} F(E_{1j}
{\rm tr}_n(X^*X)E_{j1})\right)^{1/2} \label{eq3.12}\\
&&\times\ \left(\sum^n_{j=1} g(E_{1j} \
tr_n(YY^*)E_{j1}) + G(Y^*Y)\right)^{1/2},\nonumber
\end{eqnarray}
using normality of $F$ and $g$. Since $\eta$ and $\nu$ were arbitrary,
(\ref{eq3.12}) immediately implies that
\begin{eqnarray} \qquad\quad
\|\psi(X,Y)\| &\le& C(\|XX^*\| + n {\rm tr}_n(X^*X))^{1/2} (n {\rm tr}_n(YY^*) + \|Y^*Y\|)^{1/2}\label{eq3.13}
\\
&=&  C\rho_n(X) \rho_n(Y),\nonumber
\end{eqnarray}
completing the proof.
\enddemo

\numbereddemo{{R}emark}\label{rem3.2}
The inequality (\ref{eq3.12}) implies that
\begin{equation}\label{eq3.14}
|\langle\psi(X,Y)\eta,\nu\rangle| \le C(f(XX^*) + n {\rm tr}_n(X^*X))^{1/2} (G(Y^*Y) + n {\rm tr}_n(YY^*))^{1/2}
\end{equation} 
for $X,Y\in \bb M_n(\cl M)$, which is exactly of Grothendieck type. The normal
states $F$ and $g$ have both been replaced by $n {\rm tr}_n$. The
type of averaging argument employed above may be found in \pagebreak \cite{E}.   \enddemo 

We now come to the main result of this section, a multilinear inequality which
builds on the bilinear case of Lemma~\ref{lem3.1}. We will use three versions
$\{\psi_i\}_{i=1}^3$ of the map $\psi$ in the previous lemma, with various
values of the constant~$C$. The multilinearity of $\phi$ below will guarantee
that each map satisfies the first hypothesis of Lemma \ref{lem3.1}.

\proclaim{Theorem}\label{thm3.3}
Let $\cl M\subseteq B(H)$ be a type ${\rm II}_1$ factor and let $\cl N$ be a
hyperfinite subfactor with trivial relative commutant. If $\phi\colon \ \cl
M^k\to B(H)$ is a $k$\/{\rm -}\/linear $\cl N$\/{\rm -}\/multimodular normal map{\rm ,} then 
\begin{equation}\label{eq3.15}
\|\phi^{(n)}(X_1,\ldots, X_k)\| \le 2^{k/2}\|\phi\| \rho_n(X_1)\ldots
\rho_n(X_k)
\end{equation}
for all $X_1,\ldots, X_k\in \bb M_n(\cl M)$ and $n\in\bb N$.
\endproclaim

\demo{Proof}
We may assume, without loss of generality, that $\|\phi\| = 1$. We take
(\ref{eq3.1}) as our starting point, and we will deal with the outer and inner
variables separately. Define, for $X,Y\in \bb M_n(\cl M)$,
\begin{eqnarray}\label{eq3.16}
\psi_1(X,Y)& = &\phi^{(n)}(X^*E_{11}, E_{11}X_2E_{11},\ldots, E_{11}X_kE_{11})^*\\
&&\times\ 
\phi^{(n)}(YE_{11}, E_{11}X_2E_{11},\ldots, E_{11}X_kE_{11}), \nonumber
\end{eqnarray}
where we regard $X_2,\ldots, X_k \in \bb M_n(\cl M)$ as fixed. Then
(\ref{eq3.1}) implies that
\begin{equation}\label{eq3.17}
\|\psi_1(XE_{11},E_{11}Y)\| \le \|X_2\|^2 \ldots \|X_k\|^2 \|X\| \|Y\|,
\end{equation}
and (\ref{eq3.2}) is satisfied. Taking $C$ to be $\|X_2\|^2\ldots \|X_k\|^2$
in Lemma~\ref{lem3.1} we see that
\begin{eqnarray}
\|\psi_1(X,Y)\| &=& \|\psi_1(E_{11}X,YE_{11})\|\label{eq3.18}
\\
&\le& \|X_2\|^2\ldots \|X_k\|^2 \rho_n(E_{11}X)\rho_n(YE_{11}).
\nonumber\end{eqnarray}
Now
\begin{eqnarray}
\rho_n(E_{11}X) &=& (\|E_{11}XX^*E_{11}\| + n {\rm 
tr}_n(E_{11}XX^*E_{11}))^{1/2}\label{eq3.19}
\\
&\le& 2^{1/2}\|X\|,\nonumber
\end{eqnarray}
since $\tr_n(E_{11}) = n^{-1}$, and a similar estimate holds for
$\rho_n(YE_{11})$. If we replace $X$ by $X^*_1$ and $Y$ by $X_1$ in
(\ref{eq3.18}), then (\ref{eq3.16}) and (\ref{eq3.19}) combine to give
\begin{equation}\label{eq3.20}
\|\phi^{(n)}(X_1E_{11}, E_{11}X_2E_{11},\ldots, E_{11}X_kE_{11})\|\le 2^{1/2}
\|X_1\| \|X_2\| \ldots \|X_k\|.\enspace
\end{equation}

Now consider the bilinear map
\begin{equation}\label{eq3.21}
\psi_2(X,Y) = \phi^{(n)}(X, YE_{11}, E_{11}X_3E_{11},\ldots, E_{11}X_kE_{11})
\end{equation}
where $X_3,\ldots X_k$ are fixed. By (\ref{eq3.20}), this map satisfies
(\ref{eq3.3}) with $C=2^{1/2}\|X_3\|\ldots\|X_k\|$, and multimodularity of
$\phi$ ensures that (\ref{eq3.2}) holds. By Lemma~\ref{lem3.1},
\begin{eqnarray}
\|\psi_2(X,Y)\| &= &\|\psi_2(X,YE_{11})\|\label{eq3.22}
\\
&\le &2^{1/2}\|X_3\| \ldots \|X_k\|\rho_n(X) \rho_n(YE_{11})\nonumber\\
&\le &2\|X_3\| \ldots \|X_k\|\rho_n(X)\|Y\|.\nonumber
\end{eqnarray}
Replace $X$ by $X_1$ and $Y$ by $X_2$ to obtain
\begin{equation}\label{eq3.23}
\|\phi^{(n)}(X_1,X_2E_{11},E_{11}X_3E_{11},\ldots, E_{11}X_kE_{11})\| \le
2\rho_n(X_1)\|X_2\| \ldots \|X_k\|.
\end{equation}
We repeat this step $k-2$ times across each succeeding consecutive pair of
variables, gaining a factor of $2^{1/2}$ each time and replacing each
$\|X_i||$ by $\rho_n(X_i)$, until we reach the inequality
\begin{equation}\label{eq3.24}
\|\phi^{(n)}(X_1,X_2,\ldots, X_{k-1},X_kE_{11})\| \le 2^{k/2}
\rho_n(X_1)\ldots \rho_n(X_{k-1})\|X_k\|.\quad
\end{equation}

To complete the proof, we now define
\begin{equation}\label{eq3.25}
\psi_3(X,Y) = \phi^{(n)}(X_1,\ldots, X_{k-1},X) \phi^{(n)}(X_1,\ldots, X_{k-1},
Y^*)^*,\hskip.45in
\end{equation}
where $X_1,\ldots, X_{k-1}$ are fixed. We may apply Lemma~\ref{lem3.1} with
\vglue3pt
\centerline{$C=2^k \rho_n(X_1)^2  \ldots \rho_n(X_{k-1})^2$}

\noindent  to obtain
\begin{equation}\label{eq3.26}
\|\psi_3(X,Y)\| \le C\rho_n(X) \rho_n(Y).
\end{equation}
Put $X = X_k$ and $Y = X^*_k$. Then (\ref{eq3.25}) and (\ref{eq3.26}) give the
estimate 
\begin{equation}\label{eq3.27}
\|\phi^{(n)}(X_1,\ldots, X_k)\| \le 2^{k/2} \rho_n(X_1) \ldots \rho_n(X_k),
\end{equation}
as required, since $\rho_n(X^*_k) = \rho_n(X_k)$.
\enddemo

We will use Theorem \ref{thm3.3} subsequently in a modified form which we now
state.

\vglue-18pt

\phantom{hi}

\proclaim{{C}orollary}\label{cor3.4}
Let $\cl M\subseteq B(H)$ be a type ${\rm II}_1$ factor and let $\cl N$ be a
hyperfinite subfactor with trivial relative commutant. Let $n\in \bb N${\rm ,} let
$P\in \bb M_n(\cl M)$ be a projection of trace $n^{-1}${\rm ,} and let $\phi\colon \
\cl M^k\to B(H)$ be a $k$\/{\rm -}\/linear $\cl N$\/{\rm -}\/multilinear normal map. Then{\rm ,} for
$X_1,\ldots, X_k\in \bb M_n(\cl M)${\rm ,}
\begin{equation}\label{eq3.28}
\|\phi^{(n)}(X_1P,\ldots, X_kP)\| \le 2^k\|\phi\| \|X_1\|\ldots \|X_k\|.
\end{equation}
\endproclaim

{\it Proof}.
For $1\le i\le k$,
\begin{eqnarray}
\rho_n(X_iP) &=& (\|PX^*_iX_iP\| + n {\rm tr}_n(PX^*_iX_iP))^{1/2}
\label{eq3.29}
\\ &\le& (\|X^*_iX_i\| (1+ n {\rm tr}_n(P)))^{1/2}\nonumber\\
&=& 2^{1/2} \|X_i\|.\nonumber
\end{eqnarray}
The result follows immediately from (\ref{eq3.15}) with each $X_i$ replaced by
$X_iP$.
\hfill\qed
 
\vglue-6pt
\section{Joint continuity in the ${\|\cdot\|_2}$-norm}

 \vglue-2pt

There is an extensive literature on the topic of joint and separate 
continuity of functions of two variables (see \cite{BN}, \cite{N} and the references
therein) with generalizations to the multivariable case. In this section we
consider an $n$-linear map $\phi\colon \ \cl M\times \ldots\times \cl M\to \cl
M$ on a type ${\rm II}_1$ factor $\cl M$ which is ultraweakly continuous (or normal)
separately in each variable. The restriction of $\phi$ to the closed unit ball
will be shown to be separately continuous when both range and domain have the
$\|\cdot\|_2$-norm, and from this we will deduce joint continuity in the same
metric topology. Many such joint continuity results hinge on the Baire category
theorem, and this is true of the following lemma, which we quote as a special
case of a result from \cite{BN}, and which also can be found in
\cite[p.\ 163]{Ro}. Such theorems stem from \cite{B}. 

\proclaim{Lemma}\label{lem4.1}
Let $\cl X,\cl Y$ and $\cl Z$ be complete metric spaces{\rm ,} and let $f\colon \
\cl X\times \cl Y\to \cl Z$ be continuous in each variable separately. For
each $y_0\in \cl Y${\rm ,} there exists an $x_0\in \cl X$ such that $f(x,y)$ is
jointly continuous at $(x_0,y_0)$.
\endproclaim

We now use this lemma to obtain a joint continuity result which is the first
step in an induction argument. Let $B$ denote the closed unit ball of a type
${\rm II}_1$ factor $\cl M$, to which we give the metric induced by the
$\|\cdot\|_2$-norm. Then $B$ is a complete metric space. We assume that
multilinear maps $\phi$ below satisfy $\|\phi\| \le 1$, so that $\phi$ maps
$B\times\ldots\times B$ into $B$. The $k^{\rm th}$ copy of $B$ in such a
Cartesian product will be written as $B_k$.

\proclaim{Lemma}\label{lem4.2}
Let $\phi\colon \ \cl M\times \cl M\to \cl M${\rm ,} $\|\phi\|\le 1${\rm ,} be a bilinear
map which is separately continuous in the $\|\cdot\|_2$-norm on $B_1\times
B_2$. Then $\phi\colon \ B_1\times B_2\to  B$ is jointly continuous in the
$\|\cdot\|_2$\/{\rm -}\/norm.
\endproclaim

\demo{Proof}
If we apply Lemma \ref{lem4.1} with $y_0$ taken to be 0, then there exists
$a\in B$ such that the restriction of $\phi$ to $B_1\times B_2$ (which we also
write as $\phi$) is jointly continuous at $(a,0)$. We now prove joint
continuity at (0,0), first under the assumption that $a\ge 0$, and then
deducing the general case from this. Suppose, then, that $a\ge 0$.

Consider sequences $\{h_n\}^\infty_{n=1}\in B_1$ and $\{k_n\}^\infty_{n=1} \in
B_2$, both having limit 0 in the $\|\cdot\|_2$-norm. If $h_n\ge 0$, then
$a-h_n \in B_1$ since for positive elements
\begin{equation}\label{eq4.1}
\|a-h_n\| \le \max\{\|a\|, \|h_n\|\} \le 1.
\end{equation}
Thus $\{(a-h_n,k_n)\}^\infty_{n=1}$ converges to $(a,0)$ in $B_1\times B_2$.
Since
\begin{eqnarray}
\|\phi(h_n,k_n)\|_2 &=&\|\phi((h_n-a) + a,k_n)\|_2\label{eq4.2}
\\
&\le& \|\phi(a-h_n,k_n)\|_2 + \|\phi(a,k_n)\|_2,\nonumber
\end{eqnarray}
we see that
\begin{equation}\label{eq4.3}
\lim_{n\to\infty}\|\phi(h_n,k_n)\|_2 = 0
\end{equation}
from joint continuity \pagebreak at $(a,0)$.

Now suppose that each $h_n$ is self-adjoint, and write $h_n = h^+_n-h^-_n$
with $h^+_nh^-_n = 0$, and $h^\pm_n\ge 0$. Then $\|h_n\|^2_2 = \|h^+_n\|^2_2 +
\|h^-_n\|^2_2$, so $h^\pm_h\in B_1$ and $\lim\limits_{n\to\infty}
\|h^\pm_n\|_2  = 0$. This shows that
\begin{equation}\label{eq4.4}
\lim_{n\to\infty} \phi(h_n,k_n) = \lim_{n\to\infty} \phi(h^+_n,k_n) -
\lim_{n\to\infty} \phi(h^-_n,k_n) = 0
\end{equation}
for a self-adjoint sequence in the first variable. This easily extends to a
general sequence from $B_1$ by taking real and imaginary parts. Thus $\phi$
is jointly continuous at (0,0) when $a\ge 0$.

For the general case, take the polar decomposition $a=bu$ with $b\ge 0$ and
$u$ unitary, which is possible because $\cl M$ is type ${\rm II}_1$. Then the map
$\psi(x,y) = \phi(xu,y)$ is jointly continuous at $(b,0)$, and thus at (0,0)
from above. Since $\phi(x,y) = \psi(xu^*,y)$, joint continuity of $\phi$ at
(0,0) follows immediately.

We now show joint continuity at a general point $(a,b)\in B_1\times B_2$. If
$\lim\limits_{n\to\infty} (a_n,b_n) = (a,b)$ for a sequence in $B_1\times
B_2$, then the equations
\begin{equation}\label{eq4.5}
a_n = a+2h_n,\quad b_n = b+2k_n
\end{equation}
define a sequence $\{(h_n,k_n)\}^\infty_{n=1}$ in $B_1\times B_2$ convergent
to (0,0). Then
\begin{equation}\label{eq4.6}
\phi(a_n,b_n) - \phi(a,b) = 2\phi(h_n,b) + 2\phi(a,k_n) + 4\phi(h_n,k_n),
\end{equation}
and the right-hand side converges to 0 by  joint continuity at (0,0) and
separate continuity in each variable. This shows joint continuity at $(a,b)$.
\enddemo

\proclaim{Proposition}\label{pro4.3}
Let $\phi\colon \ \cl M\times \ldots\times \cl M\to \cl M${\rm ,} $\|\phi\|\le 1${\rm ,}
be a bounded $n$\/{\rm -}\/linear map which is separately continuous in the
$\|\cdot\|_2$\/{\rm -}\/norm on $B_1\times\ldots\times B_n$. Then $\phi\colon \
B_1\times\ldots\times B_n\to B$ is jointly continuous in the
$\|\cdot\|_2$-norm.
\endproclaim

\demo{Proof}
The case $n=2$ is Lemma \ref{lem4.2}, so we proceed inductively and assume
that the result is true for all $k\le n-1$. Then consider a separately
continuous $\phi\colon \ B_1\times\ldots\times B_n\to B$. If we fix the first
variable then the resulting $(n-1$)-linear map is jointly continuous on $B_2
\times\ldots\times B_n$ by the induction hypothesis. If we view this Cartesian
product as $B_1\times (B_2\times \ldots\times B_n)$, then we have separate
continuity, so Lemma~\ref{lem4.1} ensures that there exists $a\in B_1$ so that
$\phi$ is jointly continuous at $(a,0,\ldots, 0)$. We may then follow the
proof of Lemma~\ref{lem4.2} to show firstly that $\phi$ is jointly continuous
at $(0,\ldots, 0)$, and subsequently that $\phi$ is jointly continuous at a
general point $(a_1,\ldots, a_n)$, using the induction hypothesis.
\enddemo

\proclaim{Theorem}\label{thm4.4}
Let $\phi\colon \ \cl M\times\ldots\times \cl M\to \cl M${\rm ,} $\|\phi\|\le 1${\rm ,} be
separately normal in each variable. Then the restriction of $\phi$ to $B_1
\times\ldots\times B_n$ is jointly continuous in the \pagebreak $\|\cdot\|_2$\/{\rm -}\/norm.
\endproclaim

\demo{Proof}
If we can show that the restriction of $\phi$ is separately continuous in the
$\|\cdot\|_2$-norm, then the result will follow from Proposition~\ref{pro4.3}.
By fixing all but one of the variables, we reduce to the case of a normal map
$\psi\colon \ \cl M\to \cl M$, $\|\psi\|\le 1$. By \cite[5.4.3]{SS1}, there
exist normal states $f,g\in \cl M_*$ such that
\begin{equation}\label{eq4.7}
\|\psi(x)\|_2 \le f(x^*x)^{1/2} + g(xx^*)^{1/2},\qquad x\in\cl M.
\end{equation}
We may suppose that $\cl M$ is represented in standard form, so that every
normal state is a vector state. Thus choose $\xi,\eta\in L^2(\cl M, \tr)$ such
 that 
\begin{equation}\label{eq4.8}
f(x^*x) = \langle x^*x\xi, \xi\rangle = \|x\xi\|^2_2
\end{equation}
and
\begin{equation}\label{eq4.9}
g(xx^*) = \langle xx^*\eta,\eta\rangle = \|x^*\eta\|^2_2. 
\end{equation}

Consider now a sequence $\{x_n\}^\infty_{n=1}\in B$ which converges to $x\in
B$ in the $\|\cdot\|_2$-norm. Given $\vp>0$, choose $y,z\in \cl M$ such that
\begin{equation}\label{eq4.10}
\|\xi-y\|_2, \quad \|\eta-z\|_2<\vp.
\end{equation}
Then (\ref{eq4.7})--(\ref{eq4.10}) combine to give
\begin{eqnarray}
\|\psi(x-x_n)\|_2 &\le& \|(x-x_n)\xi\|_2 + \|(x-x_n)^*\eta\|_2\label{eq4.11}
\\
&\le &\|(x-x_n)y\|_2 + \|(x-x_n)^*z\|_2 + 4\vp\nonumber\\
&\le& \|y\| \|x-x_n\|_2 + \|z\| \|x-x_n\|_2 + 4\vp.\nonumber
\end{eqnarray}
Thus, from (\ref{eq4.11}),
\begin{equation}\label{eq4.12}
\limsup_{n\ge 1}\|\psi(x-x_n)\|_2 \le 4\vp,
\end{equation}
and since $\vp>0$ was arbitrary, we conclude that $\lim\limits_{n\to \infty}
\|\psi(x-x_n)\|_2  = 0$. This proves the result.
\enddemo

\proclaim{{C}orollary}\label{cor4.5}
Let ${\cl M}\subseteq B(H)$ be a type ${\rm II}_1$ factor and let 
$\phi \colon {\cl M}\times\ldots \times {\cl M} \to B(H)${\rm ,} $\|\phi \| \leq 1${\rm ,}
be a bounded $n$\/{\rm -}\/linear map which is separately normal in each variable. For
an arbitrary pair of unit vectors $\xi,\,\eta \in H${\rm ,} 
the $n$\/{\rm -}\/linear form \begin{equation}\label{eq4.12a}
\psi(x_1,\ldots,x_n)=\langle \phi(x_1,\ldots,x_n)\xi,\eta \rangle,\qquad x_i
\in \cl M,
\end{equation}
is jointly continuous in the $\|\cdot\|_2$\/{\rm -}\/norm when restricted to 
$B_1\times\ldots\times B_n$.
\endproclaim

\demo{Proof}
View $\psi$ as having range in ${\bb C}1 \subseteq \cl M$, and apply Theorem
\ref{thm4.4}.
\enddemo

\numbereddemo{{R}emark}\label{rem4.5}
Restriction to the unit ball is necessary in the previous results. The
bilinear map $\phi(x,y) = xy$, $x,y\in \cl M$, is separately normal, but if we
take a sequence of projections $p_n\in \cl M$ of trace $n^{-4}$, then
$\lim\limits_{n\to\infty} \|np_n\|_2 = 0$, but $\|\phi(np_n, np_n)\|_2 =
n^2({\rm tr}(p_n))^{1/2} = 1$. This shows that $\phi$ is not jointly continuous in
the $\|\cdot\|_2$-norm for the whole of $\cl M$. However, a simple scaling
argument shows that $\phi$ may have arbitrary norm and that restriction to the
closed ball of any finite radius allows the same conclusion concerning joint
continuity.

If $\cl M$ is faithfully represented on a Hilbert space $H$, then the
ultrastrong$^*$ topology is defined by the family of seminorms
\begin{equation}\label{eq4.13}
x\mapsto \left(\sum^\infty_{n=1}\|x\xi_n\|^2 + \|x^*\xi_n\|^2\right)^{1/2},
\quad \xi_n\in H, \quad \sum^\infty_{n=1}\|\xi_n\|^2 < \infty,\qquad x\in\cl M.
\end{equation}
Thus convergence of a net $\{x_\lambda\}_{\lambda\in\Lambda}$ to $x$ in the
ultrastrong$^*$ topology is equivalent to ultraweak convergence of the nets
$$\{(x-x_\lambda)(x-x_\lambda)^*\}_{\lambda\in\Lambda}\quad \hbox{and}\quad
\{(x-x_\lambda)^*(x-x_\lambda)\}_{\lambda\in\Lambda}$$
to 0, showing that the ultrastrong$^*$ topology is independent of the
particular representation. By \cite[III.5.3]{T}, this topology, when
restricted to the unit ball of $\cl M$, equals the topology arising from the
$\|\cdot\|_2$-norm. Thus the conclusion of Theorem~\ref{thm4.4} could have
been stated  as the joint ultrastrong$^*$ continuity of $\phi$ when restricted
to closed balls of finite radius. In \cite{A}, Akemann proved
the equivalence of continuity in the ultraweak and ultrastrong$^*$
topologies for bounded maps restricted to balls, so these results
give another proof of Theorem \ref{thm4.4}. We have preferred to
argue directly from Grothendieck's inequality.  \enddemo

\section{Hyperfinite subfactors}\label{sec5}

In \cite{Po}, Popa showed the existence of a hyperfinite subfactor $\cl N$ of a
separable factor $\cl M$ with trivial relative commutant $(\cl N'\cap \cl M = \bb C1)$.
In this section we use Popa's result to construct such a subfactor with some
additional properties in the case that $\cl M$ has property $\Gamma$. The
second lemma below is part of the inductive step in the main theorem. We begin
with a technical result which is a special case of a more general result in \cite[Prop. 1.11]{PP}.
In our situation the proof is short and so we include it for completeness. 

\proclaim{Lemma}\label{lem5.1}
Let $\cl M$ and $\cl N$ be type ${\rm II}_1$ factors and suppose that there exists a
matrix algebra $\bb M_r$ such that $\cl M$ is isomorphic to $\bb M_r\otimes
\cl N$. If $\cl M$ has property $\Gamma${\rm ,} then so too does $\cl N$.
\endproclaim

\demo{Proof}
Fix a free ultrafilter $\omega$ on $\bb N$, and let $\cl M^\omega$ denote the
resulting ultraproduct factor, which contains a naturally embedded copy of
$\cl M$ with relative commutant denoted $\cl M_\omega$. Then $\cl M$ has
property $\Gamma$ if and only if $\cl M_\omega \ne \bb C1$ (\cite{Co}). Since
$\cl M^\omega$ is isomorphic to $\bb M_r\otimes \cl N^\omega$, and $\cl
M_\omega$ is then isomorphic to $I_r\otimes \cl N_\omega$, the result follows.
\enddemo

\proclaim{Lemma}\label{lem5.2}
Let $\cl M$ be type ${\rm II}_1$ factor with property $\Gamma$ and let $\cl M = \bb
M_r\otimes \cl N$ be a tensor product decomposition of $\cl M$. Given
$x_1,\ldots, x_k\in \cl M${\rm ,} $n\in\bb N${\rm ,} and $\vp>0${\rm ,} there exists a set of
orthogonal projections $\{p_i\}^n_{i=1}\in \cl N${\rm ,} each of trace $n^{-1}${\rm ,}
such that
\begin{equation}\label{eq5.1}
\|[1\otimes p_i, x_j]\|_2 < \vp, \qquad 1\le i \le n,\quad 1\le j \le k.
\end{equation}
\endproclaim

\demo{Proof}
Write each $x_j$ as an $r\times r$ matrix over $\cl N$, and let
$\{y_i\}^{kr^2}_{i=1}$ be a listing of all the resulting matrix entries. 
By Lemma \ref{lem5.1}, $\cl N$ has property $\Gamma$, so given
$\delta>0$ we can find a set $\{p_i\}^n_{i=1}\in\cl N$ of orthogonal
projections of trace $n^{-1}$ satisfying
\begin{equation}\label{eq5.2}
\|[p_i,y_j]\|_2 <\delta,\qquad 1\le i \le n,\quad 1\le j \le kr^2,
\end{equation}
from \cite{D}. It is clear that (\ref{eq5.1}) will hold for $\delta <
r^{-2}\vp$.
\enddemo

Since the projections found above have equal trace, they may be viewed as the
minimal projections on the diagonal of an $n\times n$ matrix subalgebra of
$\cl N$; we will use this subsequently.

\proclaim{Theorem}\label{thm5.3}
Let $\cl M$ be a type ${\rm II}_1$ factor with separable predual and with property
$\Gamma$. Then there exists a hyperfinite subfactor $\cl R$ with trivial
relative commutant satisfying the following condition. Given $x_1,\ldots,
x_k\in \cl M${\rm ,} $n\in\bb N${\rm ,} and $\vp>0${\rm ,} there exist orthogonal projections
$\{p_i\}^n_{i=1}\in\cl R${\rm ,} each of trace $n^{-1}${\rm ,} such that 
\begin{equation}\label{eq5.3}
\|[p_i,x_j]\|_2 < \vp, \qquad 1\le i \le n,\quad 1\le j \le k.
\end{equation}
\endproclaim

\demo{Proof}
We will construct $\cl R$ as the ultraweak closure of an ascending union of
matrix subfactors $\cl A_n$ which we define inductively. We first fix a
sequence $\{\theta_i\}^\infty_{i=1}$ of normal states (with $\theta_1$ the
trace) which is norm dense in the set of all normal states in $\cl M_*$. We
then choose a sequence $\{m_i\}^\infty_{i=1}$ from the unit ball of $\cl M$
which is $\|\cdot\|_2$-norm dense in the unit ball. For these choices, the
induction hypothesis is 
\begin{itemize}
\item[(i)] for each $k\le n$ there exist orthogonal projections $p_1,\ldots,
p_k\in \cl A_n$, ${\rm tr}(p_i) = k^{-1}$, satisfying
\end{itemize}

\vglue-18pt

\begin{equation}\label{eq5.4}
\|[p_i,m_j]\|_2  < n^{-1},\qquad 1\le i \le k,\quad 1\le j \le n;
\end{equation}
\begin{itemize}
\item[(ii)] if $\cl U_n$ is the unitary group of $\cl A_n$ with normalized Haar
measure $du$, then \end{itemize}

\vglue-18pt

\begin{equation}\label{eq5.5}
\left|\theta_i\left(\int_{\cl U_n} um_ju^* \ du\right) - {\rm tr}(m_j)\right| <
n^{-1}
\end{equation}
\begin{itemize} \item[] for $1\le i \le n, \ 1\le j \le n$.
\end{itemize}

To begin the induction, let $\cl A_1 = \bb C1$ and let $p_1=1$, which commutes
with $m_1$, so (i) holds. The second part of the hypothesis is also satisfied
because $\theta_1$ is the trace. Now suppose that $\cl A_{n-1}$ has been
constructed. We apply Lemma~\ref{lem5.2} $n$ times to the set $\{m_1,\ldots,
m_n\}$, taking $\vp$ to be $n^{-1}$ and $k$ to be successively $1,2,\ldots,
n$. At the $k^{\rm th}$ step we acquire a copy of $\bb M_k$, leading to a
matrix algebra  
\begin{equation}\label{eq5.5a}
\cl B_n = \cl A_{n-1} \otimes \bb M_{1} \otimes
\bb M_{2}\otimes  \cdots\otimes \bb M_{n}=\cl A_{n-1} \otimes \bb M_{n!}
\end{equation} 
containing sets of projections which satisfy (i).

Now decompose $\cl M$ as $\cl B_n\otimes\cl N$ for some type ${\rm II}_1$ factor $\cl
N$, and choose a hyperfinite subfactor $\cl S\subseteq \cl N$ with trivial
relative commutant, \cite{Po}. There exists an ascending sequence $\{\cl
F_r\}^\infty_{r=1}$ of matrix subalgebras of $\cl S$ whose union is
ultraweakly dense in $\cl S$, as is $\bigcup\limits_{r\ge 1}\cl B_n\otimes \cl
F_r$ in $\cl B_n \otimes \cl S$. Let $\cl V_r$ denote the unitary group of
$\cl B_n \otimes \cl F_r$ with normalized Haar measure $dv$. Since $\cl
B_n\otimes \cl S$ has trivial relative commutant in ${\cl M}$, a standard
computation (see, for example, \cite[5.4.4]{SS1}) shows that
\begin{equation}\label{eq5.6}
\lim_{r\to\infty} \int_{\cl V_r} vxv^* \ dv = {\rm tr}(x)1
\end{equation}
ultraweakly for all $x\in\cl M$. Since each $\theta_i$ is normal, we may select
$r$ so large that \begin{equation}\label{eq5.7}
\left|\theta_i\left(\int_{\cl V_r} vm_jv^* \ dv\right) - {\rm tr}(m_j)\right| <
n^{-1} 
\end{equation}
for $1\le i \le n$ and $1\le j \le n$. For this choice of $r$, define $\cl
A_n$ to be $\cl B_n\otimes \cl F_r$. Now both (i) and (ii) are satisfied.

Let $\cl R\subseteq \cl M$ be the ultraweak closure of the union of the $\cl
A_n$'s. We now verify that (\ref{eq5.3}) holds for a given set
$\{x_1,\ldots,x_k\} \in \cl M$, $n \in \bb N$ and $\vp > 0$. Let 
$\delta = \vp /3$ and, without loss of generality, assume that $\|x_j\| \leq
1$ for $1 \leq j \leq k$. Then choose elements $m_{n_j}$ from the sequence so
that
\begin{equation}\label{eq5.7a}
\|x_j - m_{n_j}\|_2 < \delta, \qquad 1 \leq j \leq k.
\end{equation}
Now select $r \in \bb N$ to be so large that
\begin{equation}\label{eq5.7b}
r > \delta^{-1},\ n,\ {\rm max}\,\{n_j:\,1 \leq j \leq k\}.
\end{equation}
By (i) (with $n$ and $r$ replacing respectively $k$ and $n$), there exist
orthogonal
projections $p_i \in {\cl A}_r \subseteq \cl R$, $1 \leq i \leq n$, each of
trace $n^{-1}$, such that 
\begin{equation}\label{eq5.7c}
\|[p_i,m_{n_j}]\|_2 < r^{-1} < \delta, \qquad 1\le i \le n,\quad 1\le j \le k. \pagebreak
\end{equation}
Then (\ref{eq5.7a}), (\ref{eq5.7c}) and the triangle inequality give
\begin{equation}\label{eq5.7d}
\|[p_i,x_j]\|_2 < 3\delta = \vp, \qquad 1\le i \le n,\quad 1\le j \le k,
\end{equation}
as required.
It remains to show that
$\cl R' \cap \cl M = \bb C1$, which will also show that ${\mathcal R}$ is a factor.

Consider $x\in \cl R' \cap \cl M$, which we may assume to be of unit norm.
Then choose a subsequence $\{m_{n_j}\}^\infty_{j=1}$ converging to $x$ in the
$\|\cdot\|_2$-norm. We note that $\|x-m_{n_j}\|\le 2$, so this sequence
converges to $x$ ultraweakly, and $\lim\limits_j {\rm tr}(m_{n_j}) = {\rm tr}(x)$, by
normality of the trace. The inequality
\begin{eqnarray}
\left\|\left(\int_{\cl U_{n_j}} um_{n_j}u^* \ du\right)-x\right\|_2 &=&
\left\|\int_{\cl U_{n_j}} u(m_{n_j}-x)u^* \ du\right\|_2\label{eq5.8}
\\
&\le& \|m_{n_j}-x\|_2,\nonumber
\end{eqnarray}
which is valid because $x\in \cl R'$, shows that these integrals also converge
ultraweakly to $x$. For any fixed value of $i$, the sequence
$\left\{\theta_i\left(\int_{\cl U_{n_j}} um_{n_j} u^* \
du\right)\right\}^\infty_{j=1}$\break
\vglue-3pt\noindent converges to $\theta_i(x)$, since each
$\theta_i$ is ultraweakly continuous, and also to ${\rm tr}(x)$, by (\ref{eq5.7}).
This shows that $\theta_i(x)={\rm tr}(x)$ for each $i \geq 1$. By norm density of
$\{\theta_i\}^\infty_{i=1}$ in the set of normal states, we conclude that $x =
{\rm tr}(x)1$, and so $\cl R$ has trivial relative commutant in $\cl M$.
\enddemo 

\numbereddemo{{R}emark}\label{rem5.4}
We note, from the construction of the ${\cl A}_n$'s, that the
projections in the previous theorem are contained in a Cartan masa in $\cl R$.
It is not clear whether this is a masa in $\cl M$ in general (and we would not
expect it to be Cartan in $\cl M$). We do not pursue this point as it will not
be needed subsequently.   \enddemo 

\section{The separable predual case}\label{sec6}

In this section we show that the cohomology groups $H^k(\cl M,\cl M)$
and $H^k(\cl M,B(H))$, $k\ge
2$, are 0 for any type ${\rm II}_1$ factor ${\cl M}\subseteq B(H)$ with property
$\Gamma$ and separable predual (but note that we place no restriction on $H$).
The general case is postponed to the next section. We will need an algebraic
lemma, for which we now establish some notation.

Let $S_k$, $k\ge 2$, be the set of nonempty subsets of $\{1,2,\ldots, k\}$,
and let $T_k$ be the collection of subsets containing $k$. The cardinalities
are respectively $2^k-1$ and $2^{k-1}$. If $\sigma\in S_k$ then we also regard
it as an element of $S_r$ for all $r>k$, and we denote its cardinality by
$|\sigma|$. We note that $S_{k+1}$ is then the disjoint union of $S_k$ and
$T_{k+1}$. If $\phi\colon \ \cl M^k\to B(H)$ is a $k$-linear map, $p\in \cl
M$ is a projection and $\sigma\in S_k$, then we define $\phi_{\sigma,p}\colon
\ \cl M^k\to B(H)$ by
\begin{equation}\label{eq6.1}
\phi_{\sigma,p}(x_1,\ldots, x_k) = \phi(y_1,\ldots, y_k),
\end{equation}
where $y_i = px_i-x_ip$ for $i\in\sigma$, and $y_i = x_i$ otherwise. For
convenience of notation we denote the commutator $[p,x_i]$ by $\hat x_i$ since
we will only be concerned with one projection at this time. For example, if
$k=3$ and $\sigma = \{2,3\}$, then
\begin{eqnarray}
\phi_{\sigma,p}(x_1,x_2,x_3) &=& \phi(x_1,px_2-x_2p,px_3-x_3p)\label{eq6.2}
\\
&=& \phi( x_1,\hat x_2,\hat x_3),\qquad x_i\in\cl M.\nonumber
\end{eqnarray}
If $\sigma \in S_k$, denote by $\ell(\sigma)$ the least integer in $\sigma$.
Then define $\phi_{\sigma,p,i}(x_1,\ldots, x_k)$ by changing the $i^{\rm th}$
variable in $\phi_{\sigma,p}$ from $x_i$ to $\hat x_i$, $1 \leq i <
\ell(\sigma)$, and replacing $\hat x_i$ by $p\hat x_i$ when $i=\ell(\sigma)$.
In the above example $\ell(\sigma)=2$, and
\begin{equation}\label{eq6.2a}
\phi_{\sigma,p,1}(x_1,x_2,x_3)=\phi(\hat x_1,\hat x_2,\hat x_3),\quad
\phi_{\sigma,p,2}(x_1,x_2,x_3)=\phi( x_1,p\hat x_2,\hat x_3).\quad
\end{equation}

\proclaim{Lemma}\label{lem6.1}
Let $p\in {\cl M} \subseteq B(H)$ be a fixed but arbitrary projection{\rm ,} and let
$\cl C_k$, $k\ge 2${\rm ,} be the class of $k$\/{\rm -}\/linear maps $\phi\colon \ \cl M^k\to
B(H)$ which satisfy
\begin{eqnarray}
\label{eq6.3}
p\phi(x_1,\ldots, x_k) &=& \phi(px_1,x_2,\ldots, x_k),\\
\label{eq6.4}
\phi(x_1,\ldots, x_ip, x_{i+1},\ldots, x_k) &=& \phi(x_1,\ldots, x_i,px_{i+1},
\ldots, x_k),
\end{eqnarray}
for $x_j\in\cl M$ and $1\le i \le k-1$. If $\phi\in\cl C_k$ then
\begin{equation}\label{eq6.5}
p\phi(x_1,\ldots, x_k) -p\phi(x_1p,\ldots, x_kp) = \sum_{\sigma\in S_k}
(-1)^{|\sigma|+1} p\phi_{\sigma,p}(x_1,\ldots, x_k).\quad
\end{equation}
Moreover{\rm ,} for each $\sigma \in S_k${\rm ,}
\begin{equation}\label{eq6.5a}
p\phi_{\sigma,p}(x_1,\ldots,x_k)=\sum^{\ell(\sigma)}_{i=1}\phi_{\sigma,p,i}
(x_1,\ldots,x_k).
\end{equation}
\endproclaim

\demo{Proof}
We will show (\ref{eq6.5}) by induction, so consider first the case $k=2$ and
take $\phi\in\cl C_2$. Then, using (\ref{eq6.3}) and (\ref{eq6.4}) 
repeatedly,
\begin{eqnarray}
\quad\qquad p\phi(x_1,x_2) &= &p\phi(px_1,x_2)\label{eq6.6}
\\
&=& p\phi(\hat x_1,x_2) + p\phi(x_1p,x_2)\nonumber\\
&=& p\phi(\hat x_1,x_2)  + p\phi(x_1p,px_2)\nonumber\\
&=& p\phi(\hat x_1,x_2) + p\phi(x_1p,\hat x_2) + p\phi(x_1p,x_2p)\nonumber\\
&=& p\phi(\hat x_1,x_2) - p\phi(\hat x_1,\hat x_2) + p\phi(px_1, \hat x_2) +
p\phi(x_1p,x_2p)\nonumber\\
&=& p\phi(\hat x_1,x_2) - p\phi(\hat x_1,\hat x_2) + p\phi(x_1,\hat x_2) +
p\phi(x_1p,x_2p)\nonumber
\end{eqnarray}
and the result follows by moving $p\phi(x_1p,x_2p)$ to the left-hand side.

Suppose now that (\ref{eq6.5}) is true for maps in $\cl C_r$ with $r<k$, and
consider $\phi\in\cl C_k$. Note that if we fix $x_k$, the resulting
map is an element of $\cl C_{k-1}$, so the induction hypothesis gives
\begin{eqnarray}\label{eq6.7}
&&\hskip-48pt p\phi(x_1,\ldots, x_k) - p\phi(x_1p,\ldots, x_{k-1}p,x_k) \\
&&\quad\qquad = \sum_{\sigma\in
S_k\backslash T_k} (-1)^{|\sigma|+1} p\phi_{\sigma,p}(x_1,\ldots, x_k). \nonumber
\end{eqnarray}
Since this is an algebraic identity, we may replace $x_k$ by $\hat x_k$ to
obtain 
\begin{eqnarray}
&&\hskip-36pt p\phi(x_1,\ldots, x_{k-1},\hat x_k) - p\phi(x_1p,\ldots, x_{k-1}p,\hat x_k)
 \label{eq6.8}\\ 
&&\qquad = \sum_{\sigma\in S_k\backslash T_k} (-1)^{|\sigma|+1} p\phi_{\sigma,p}(x_1,
\ldots, x_{k-1},\hat x_k).\nonumber
\end{eqnarray}
By (\ref{eq6.4}),
\begin{eqnarray} \qquad 
p\phi(x_1p,\ldots, x_{k-1}p,x_k) &\hskip-8pt =\hskip-8pt& p\phi(x_1p,\ldots, x_{k-1}p, px_k)
\label{eq6.9}
\\
&\hskip-8pt=\hskip-8pt& p\phi(x_1p,\ldots, x_{k-1}p, \hat x_k) + p\phi(x_1p,\ldots, x_kp).\nonumber
\end{eqnarray}
Now use (\ref{eq6.8}) to replace $p\phi(x_1p,\ldots, x_{k-1}p, \hat x_k)$ in
(\ref{eq6.9}), and add the resulting equation to (\ref{eq6.7}). After
rearranging, we obtain
\begin{eqnarray}
 \label{eq6.10}
&&\hskip-.5in p\phi(x_1,\ldots, x_k) - p\phi(x_1p,\ldots, x_kp)\\[5pt]
& & = \ \sum_{\sigma\in
S_k\backslash T_k} (-1)^{|\sigma|+1} p\phi_{\sigma,p}(x_1,\ldots,x_k)
\nonumber\\
&  &\qquad -\ \sum_{\sigma\in S_k\backslash T_k} (-1)^{|\sigma|+1} p\phi_{\sigma,p}
(x_1,\ldots, x_{k-1}, \hat x_k)\nonumber\\
& &\qquad +\ p\phi(x_1,\ldots,x _{k-1},\hat x_k),\nonumber
\end{eqnarray}
and the right-hand side of (\ref{eq6.10}) is equal to $\sum\limits_{\sigma\in
S_k} (-1)^{|\sigma|+1} p\phi_{\sigma,p}(x_1,\ldots, x_k)$. This completes the
inductive step.

We now prove the second assertion. The idea is to bring the projection in on
the left, then past each variable (introducing a commutator each time) until
the first existing commutator is reached. To avoid technicalities we
illustrate this in the particular case of $k=3$ and $\sigma=\{2,3\}$. We use
(\ref{eq6.3}) and (\ref{eq6.4}) to move $p$ to the right, and the general
procedure should then be clear. Thus
\begin{eqnarray}
p\phi(x_1,\hat x_2,\hat x_3) &= &\phi(px_1,\hat x_2,\hat x_3)\label{eq6.10a}\\
&=&\phi(\hat x_1,\hat x_2,\hat x_3)+\phi(x_1p,\hat x_2,\hat x_3)\nonumber\\
&=&\phi(\hat x_1,\hat x_2,\hat x_3)+\phi(x_1,p\hat x_2,\hat
x_3),\nonumber
\end{eqnarray}
as required.
\enddemo

\proclaim{Lemma}\label{lem6.2}
Let ${\cl M} \subseteq B(H)$ be a type ${\rm II}_1$ factor and let $\phi \colon {\cl
M}^k \to B(H)$ be a bounded $k$\/{\rm -}\/linear separately normal map. Let
$\{p_r\}_{r=1}^{\infty}$ be a sequence of projections in $\cl M$ which satisfy
{\rm (\ref{eq6.3}), (\ref{eq6.4})} and 
\begin{equation}\label{eq6.10b}
\lim_{r\to\infty} \|[p_{r},x]\|_2 = 0,\qquad x \in \cl M.
\end{equation}
Then for each $\sigma \in S_k${\rm ,} each integer $i \leq \ell(\sigma)$ and each
pair of unit vectors $\xi,\,\eta \in H${\rm ,}
\begin{equation}\label{eq6.10c}
\lim_{r\to\infty}\langle \phi_{\sigma,p_r,i}(x_1,\ldots,x_k)\xi,\eta\rangle=0.
\end{equation}
\endproclaim

\demo{Proof}
By Lemma \ref{lem6.1}, each variable in $\phi_{\sigma,p_r,i}$ is one of three
types, and at least one of the latter two must occur: $x_j$, $p_rx_j-x_jp_r$
and $p_r(p_rx_j-x_jp_r)$. Thus, as $r \to \infty$, the variables either remain
the same (first type)
or tend to 0 in the $\|\cdot\|_2$-norm (second and third types), by
hypothesis. The result follows from the joint continuity of Corollary
\ref{cor4.5}. \enddemo

We now come to the main result of this
section, the vanishing of cohomology for property $\Gamma$ factors with
separable predual. The heart of the proof is to show complete boundedness of
certain multilinear maps and we state this as a separate theorem.

\proclaim{Theorem}\label{thm6.2a}
Let ${\cl M}\subseteq B(H)$ be a type ${\rm II}_1$ factor with property $\Gamma$ and
a separable predual. Let ${\cl R} \subseteq {\cl M}$ be a hyperfinite
subfactor with trivial relative commutant and satisfying the conclusion of
Theorem {\rm \ref{thm5.3}.} Then a bounded
$k$\/{\rm -}\/linear $\cl R$\/{\rm -}\/multimodular separately normal map $\phi\colon {\cl M}^k
\to B(H)$ is completely bounded and $\|\phi\|_{\rm cb} \leq 2^k\,\|\phi\|$.
\endproclaim

\demo{Proof}
Fix an integer $n$, and a set $X_1,\ldots,
X_k \in \bb M_n(\cl M)$. By Theorem~\ref{thm5.3}, we may find sets of
orthogonal projections $\{p_{i,r}\}^n_{i=1}$, $r\ge 1$, in $\cl R$ with trace
$n^{-1}$ such that for each $x\in \cl M$
\begin{equation}\label{eq6.13}
\lim_{r\to\infty} \|[p_{i,r},x]\|_2 = 0,\qquad 1\le i \le n.
\end{equation}
Let $P_{i,r}\in \bb M_n(\cl M)$ be the diagonal projection $I_n\otimes
p_{i,r}$. These projections satisfy the analog of (\ref{eq6.13}) for elements
of $\bb M_n(\cl M)$.

The $n$-fold amplification $\phi^{(n)}$ of $\phi$ to $\bb M_n(\cl M)$ is an
$\bb M_n(\cl R$)-multimodular map, so (\ref{eq6.3}) and (\ref{eq6.4}) are
satisfied. Thus, for each $r\ge 1$, it follows from Lemma~\ref{lem6.1} that
\begin{eqnarray} &&\label{eq6.14}\\
&&\sum^n_{i=1} P_{i,r}\phi^{(n)}(X_1,\ldots, X_k) - \sum^n_{i=1}
\sum_{\sigma\in S_k} (-1)^{|\sigma|+1} P_{i,r} \phi^{(n)}_{\sigma, P_{i,r}}
(X_1,\ldots, X_k)\nonumber
\\[5pt]
&&\quad = \sum^n_{i=1} P_{i,r}\phi^{(n)}(X_1P_{i,r},\ldots, X_kP_{i,r})
= \sum^n_{i=1} P_{i,r}\phi^{(n)}(X_1P_{i,r},\ldots, X_kP_{i,r})P_{i,r},\nonumber
\end{eqnarray}
where the last equality results from multimodularity of $\phi^{(n)}$.
Since $\{P_{i,r}\}_{i=1}^n$ is a set of orthogonal projections for each $r
\geq 1$, the right-hand side of (\ref{eq6.14}) has norm at most
\begin{equation}\label{eq6.15} \max_{1\le i\le n}
\{\|\phi^{(n)}(X_1P_{i,r},\ldots, X_kP_{i,r})\|\} \le 2^k\,\|\phi\|\,\|X_1\|
\ldots \|X_k\|, \hskip.4in\end{equation}
using (\ref{eq3.28}) in \pagebreak Corollary \ref{cor3.4}.

Now fix an arbitrary pair of unit vectors $\xi,\,\eta \in H^n$, and apply the
vector functional $\langle \,\cdot \,\xi,\eta \rangle$ to (\ref{eq6.14}). When
we let $r \to \infty$ in the resulting equation, the terms in the double sum
tend to 0 by Lemmas \ref{lem6.1} and \ref{lem6.2}, leaving the inequality
\begin{equation}\label{eq6.15a}
|\langle  \phi^{(n)}(X_1,\ldots, X_k)\xi,\eta\rangle| \leq 
2^k\,\|\phi\|\,\|X_1\| \ldots \|X_k\|,
\end{equation}
since the projections in the first term of (\ref{eq6.14}) sum to 1.
Now $n$, $\xi$ and $\eta$ were arbitrary, so complete boundedness of $\phi$ 
follows from (\ref{eq6.15a}), as does the inequality
$\|\phi\|_{\rm cb} \leq 2^k\,\|\phi\|$.
\enddemo

In the following theorem we restrict to $k \geq 2$ since the two cases of
$k=1$ are in \cite{K,S} and \cite{C2} respectively.

 \proclaim{Theorem}\label{thm6.2}
Let ${\cl M} \subseteq B(H)$ be a type ${\rm II}_1$ factor with property $\Gamma$
and a separable predual. Then 
\begin{equation}\label{eq6.12} 
H^k(\cl M,\cl M) =
H^k(\cl M,B(H)) =
0,\qquad k\ge 2. \end{equation} \endproclaim

\demo{Proof} 
Let $\cl R$ be a hyperfinite subfactor of $\cl M$ with trivial relative
commutant and satisfying the additional property of Theorem~\ref{thm5.3}. 
We consider first the cohomology groups $H^k(\cl M,\cl M)$. By
\cite[Chapter 3]{SS1}, it suffices to consider an $\cl R$-multimodular
separately normal $k$-cocycle $\phi$, which is then
completely bounded by Theorem \ref{eq6.2a}.   It now follows from \cite{CS3}, \cite{CS4} (see also \cite[4.3.1]{SS1})
that $\phi$ is a coboundary. When $B(H)$ is the module, we appeal instead to
\cite{CES} to show that each completely bounded cocycle is a coboundary,
completing the proof.  \enddemo

\numbereddemo{{R}emark}\label{rem6.3a}
By \cite[Chapter 3]{SS1}, cohomology can be reduced to the consideration of
normal ${\cl R}$-multimodular maps which, in the case of property $\Gamma$
factors, are all completely bounded from Theorem \ref{thm6.2a}. Thus we reach
the perhaps surprising conclusion that
\begin{equation}\label{6.17a}
H^k({\cl M},{\cl X})=H^k_{\rm cb}({\cl M},{\cl X}),\qquad k \geq 1,
\end{equation}
for any property $\Gamma$ factor $\cl M$ and any ultraweakly closed 
$\cl M$-bimodule $\cl X$ lying between $\cl M$ and $B(H)$. 
\enddemo

\numbereddemo{{R}emark}\label{rem6.3}
Theorem \ref{thm6.2} shows that each normal $k$-cocycle $\phi$ may be
expressed as $\partial\psi$ where $\psi\colon \ \cl M^{k-1}\to \cl M$ (or into
$B(H)$). Lemma~3.2.4 of \cite{SS1} and the proof of Theorem~5.1 of \cite{SS2}
make it clear that $\psi$ can be chosen to satisfy
\begin{equation}\label{eq6.18}
\|\psi\| \le K_k\|\phi\|
\end{equation}
for some absolute constant $K_k$. 
\enddemo

\section{The general case}\label{sec7}

We now consider the general case of a type ${\rm II}_1$ factor $\cl M$ which has
property $\Gamma$, but is no longer required to have a separable predual. We
will, however, make use of the separable predual case of the previous section.
The connection is established by our first result.

\proclaim{Proposition}\label{pro7.1}
Let $\cl M$ be a type ${\rm II}_1$ factor with property $\Gamma${\rm ,} let $F$ be a
finite subset of $\cl M${\rm ,} and let $\phi\colon \ \cl M^k\to \cl M$ be a
bounded $k$\/{\rm -}\/linear separately normal map. Then $F$ is contained in a
subfactor $\cl M_F$ which has property $\Gamma$ and a separable predual.
Moreover{\rm ,} $\cl M_F$ may be chosen so that $\phi$ maps $(\cl M_F)^k$ into $\cl
M_F$.
\endproclaim

\demo{Proof}
We will construct inductively an ascending sequence of separable unital
$C^*$-subalgebras $\{\cl A_n\}^\infty_{n=1}$ of $\cl M$, each containing $F$,
with the following properties:
\begin{itemize}
\item[(i)] $\phi$ maps $(\cl A_n)^k$ into $\cl A_{n+1}$;

\item[(ii)] given $x_1,\ldots, x_r\in \cl A_n$ and $\vp>0$, there exists a unitary
$u\in \cl A_{n+1}$ of trace 0 such that  \end{itemize}
\vglue-22pt
\begin{equation}\label{eq7.1}
\|[x_i,u]\|_2 < \vp,\qquad 1\le i \le r;
\end{equation}
\begin{itemize}
\item[(iii)] there exists a sequence of unitaries $\{v_i\}^\infty_{i=1}$ in $\cl
A_{n+1}$ such that \end{itemize}
\vglue-22pt
\begin{equation}\label{eq7.2}
{\rm tr}(x)1 \in \ovl{\hbox{conv}}{}^{\|\cdot\|} \{v_ixv^*_i\colon \ i\ge 1\}
\end{equation}
\begin{itemize}
\item[] for each $x\in\cl A_n$.
\end{itemize}

Define $\cl A_1$ to be the separable $C^*$-algebra generated by the elements
of  $F$ and the identity element. We will only show the construction of $\cl
A_2$, since the inductive step from $\cl A_n$ to $\cl A_{n+1}$ is identical.

The restriction  of $\phi$ to $(\cl A_1)^k$ has separable range which,
together with~$\cl A_1$, generates a separable $C^*$-algebra $\cl B$. Then
$\phi$ maps $(\cl A_1)^k$ into $\cl B$. Now fix a countable sequence
$\{a_n\}^\infty_{n=1}$ which is norm dense in the unit ball of $\cl A_1$. 
For each finite subset $\sigma$ of this sequence and each integer $j$ we may
choose a trace 0 unitary $u_{\sigma,j}$ such that
\begin{equation}\label{eq7.3}
\|[a,u_{\sigma,j}]\|_2 < j^{-1},\qquad a\in\sigma.
\end{equation}
There are a countable \pagebreak number of such unitaries, so together with $\cl B$ they
generate a larger separable $C^*$-algebra $\cl C$. By the Dixmier 
approximation theorem, \cite{Dix}, we may choose a countable set of unitaries
$\{v_i\}^\infty_{i=1}\in\cl M$ so that (\ref{eq7.2}) holds  when $x$ is any
element of $\{a_n\}^\infty_{n=1}$. Then these unitaries, combined with~$\cl
C$, generate a separable $C^*$-algebra $\cl A_2$. By construction of $\cl
B$, $\phi$ maps $(\cl A_1)^k$ into~$\cl A_2$, while the second and third
properties follow from a simple approximation argument using the norm density
of $\{a_n\}^\infty_{n=1}$.

Let $\cl A_F$ be the norm closure of $\bigcup\limits_{n\ge 1} \cl A_n$, and
denote the ultraweak closure by $\cl M_F$. Then $\cl M_F$ has separable
predual  and property $\Gamma$, from (\ref{eq7.1}) and the $\|\cdot\|_2$-norm
density of $\cl A_F$ in $\cl M_F$. It remains to show that $\cl M_F$ is a
factor. If $\tau$ is a normalized normal trace on $\cl M_F$ then (\ref{eq7.2})
shows that $\tau$ and ${\rm tr}$ agree on $\cl A_F$. By normality they agree on $\cl
M_F$, so this von~Neumann algebra has a unique normalized normal trace and is
thus a factor. This completes the proof.
\enddemo

\proclaim{Theorem}\label{thm7.2}
Let ${\cl M} \subseteq B(H)$ be a type ${\rm II}_1$ factor with property $\Gamma$.
Then \begin{equation}\label{eq7.4}
H^k(\cl M,\cl M) = H^k({\cl M},B(H)) =0,\qquad k\ge 2.
\end{equation}
\endproclaim

\demo{Proof} We first consider $H^k(\cl M,\cl M)$.
By \cite[Chapter 3]{SS1}, we may restrict attention to a separately normal
$k$-cocycle $\phi$. For each finite subset $F$ of $\cl M$, let $\phi_F$ be
the restriction of $\phi$ to the subfactor $\cl M_F$ of
Proposition~\ref{pro7.1}. By Theorem~\ref{thm6.2}, there exists a
$(k-1)$-linear map $\psi_F\colon \ (\cl M_F)^{k-1}\to \cl M_F$ such that
$\phi_F = \partial\psi_F$ and there is a uniform bound on $\|\psi_F\|$
(Remark~\ref{rem6.3}). Let $\bb E_F$ be the normal conditional expectation of
$\cl M$ onto $\cl M_F$, and define $\theta_F\colon \ \cl M^{k-1}\to \cl M$ by
the composition $\psi_F\circ (\bb E_F)^{k-1}$. Any $F$ which contains a given
set $\{x_1,\ldots, x_k\}$ of elements of $\cl M$ satisfies
\begin{equation}\label{eq7.5}
\phi(x_1,\ldots, x_k) = \phi_F(x_1,\ldots, x_k) = \partial\theta_F(x_1,\ldots,
x_k).
\end{equation}
Now order the finite subsets of $\cl M$ by inclusion and take a point
ultraweakly convergent subnet of $\{\theta_F\}$ with limit $\theta\colon \ \cl
M^{k-1}\to \cl M$. It is then a simple matter to check that $\phi =
\partial\theta$, and thus $H^k(\cl M,\cl M) = 0$.

The case of $H^k({\cl M},B(H))$ is essentially the same. The only difference
is that $\psi_F$ and $\theta_F$ now map into $B(H)$ in place of ${\cl M}_F$. 
\enddemo

\numbereddemo{{R}emark}\label{rem7.3}
A more complicated  construction of $\cl M_F$ in the preceding two results
would have given the additional property that $\cl M_F\subseteq \cl M_G$
whenever $F\subseteq G$ is an inclusion of finite subsets of $\cl M$. However,
this was not needed for Theorem~\ref{thm7.2}. 
\enddemo

\end{document}